\documentclass[12pt,a4paper,oneside]{amsart}
\usepackage[latin1]{inputenc}
\usepackage{amsmath, amsfonts, amssymb,amsthm}

\author{Rafael Tesoro}
\thanks{I am grateful to Javier Cilleruelo for his suggestions and comments that improved the content of this paper.}
\address{R. Tesoro:  Departamento de Matem\'{a}ticas, Universidad
Aut\'onoma de Madrid. 28049 Madrid, Spain}
\email{rafael.tesoro@estudiante.uam.es}
\title{Sets of integers avoiding congruent subsets}
\date{24/06/2013}
\usepackage
			{hyperref}			

\usepackage{todo}				

\usepackage[a4paper]{geometry}
\usepackage{mathtools}

\newtheorem{thm}{Theorem}

\theoremstyle{definition}

\theoremstyle{remark}

\numberwithin{equation}{section}

\usepackage{dsfont}

  \newcommand\Y{\mathds{Y}}

\renewcommand{\Pr}[1]{\mathbb{P}\left (  #1 \right )}
\newcommand{\Esp}[1]{\mathbb{E}\left (  #1 \right )}

\newcommand{\ud}{\,\mathrm{d}}

\begin{document}

\begin{abstract}
We study the sets of integers $A$ that avoid any arrangement of $g$ congruent $h$-subsets (the $C_h[g]$ sets, for short), as well as the variant in which the $h$-subsets are pairwise disjoint (the weak-$C_h[g]$ sets).
For $h=g=2$ these are the Sidon sets and the weak-Sidon sets respectively.

We refine and improve upon several results by Erd\H{o}s and Harzheim.
For finite sets on the one hand we prove the upper bound $$|A| \le (g-1)^{1/h} n^{1-\frac{1}{h}} + O\left (n^{\frac{1}{2}-\frac 1{2h}}\right ),$$ for any $C_h[g]$ set $A \subset [1,n]$. On the other hand we prove that there exists a weak-$C_h[g]$ subset $A \subset [1,n]$ such that
$$|A| \gg n^{\left( 1- \frac{1}{h}\right)\left( 1- \frac{1}{g} \right)\left(1+\frac{1}{hg-1} \right)}.$$

For any infinite $C_h[g]$ sequence $A$ we prove that $$A(x) \ll x^{1-1/h} (\log x)^{-1/h}$$ for infinitely many integers $x$. 

We conclude stating several related open problems.
\end{abstract}

\maketitle

\section{Introduction}
We say that a set of integers $A$ is a $C_h[g]$ set if for any set $X$ of $h$ elements
there do not exist distinct integers $k_1,\cdots,k_g$ such that $X+k_1,\cdots ,X+k_g\subset A$.
We write \emph{weak-}$C_h[g]$ when the condition that the sets $X+k_i$ are pairwise disjoint is added.
This generalization of Sidon sets was introduced by Erd\H os and Harzheim in \cite{ErdHar-1986}.

Erd\H{o}s and Harzheim used the term \lq \lq $B_{gh}$-sequence\rq \rq.
We prefer to write $C_h[g]$ (with $C$ as in \lq \lq \textbf{c}ongruent\rq \rq) aiming to prevent confusion with the established convention of using $B_h[g]$ for a different generalization of Sidon sets.
Our notation also allows to state in a concise way some symmetries such us the fact that any $C_h[g]$ set is also a $C_g[h]$ set, and vice versa.
Indeed if $X_1, \cdots, X_g \subset A$ are congruent $h$-sets, for $\nu=1, \cdots, h$ let $Y_\nu$ be the set of the $\nu^{th}$ elements of
the sets $X_i$. Then $Y_1, \cdots, Y_h$ are congruent $g$-subsets of $A$.

In the sequel we will assume $g \ge h \ge 2.$
We will also use the following standard notation:

 $\bullet$ $f(n) \ll g(n)$ or $f(n)=O(g(n))$ means that there exists $C>0$ such that $f(n) < C g(n)$ for all $n$.

 $\bullet$ $f(x)=o(g(x))$ means that $f(x)/g(x)\to 0$ as $x\to \infty$.

 $\bullet$ $f(x) \sim g(x)$ means that $f(x)/g(x)\to 1$ as $x\to \infty$.

\medskip
Our first result is the following.
\begin{thm}  \label{thm1}
If $A \subset [1,n]$ is a $C_h[g]$ set, with $g \ge h \ge 2$,  then
\begin{equation*} \label{eq:BetterUpperBound-for-F_h(g,n)}
    |A| \le (g-1)^{1/h} n^{1-\frac{1}{h}} + O\left(n^{\frac{1}{2}-\frac 1{2h}} \right).
\end{equation*}
\end{thm}
This theorem is a refinement of the estimate $|A|\ll n^{1-1/h}$ proved by Erd\H{o}s and Harzheim.
We remark that $C_2[2]$ sets are just Sidon sets and Theorem \ref{thm1} recovers the well known upper bound
for the size of Sidon sets in $\{1,\cdots,n\}$ obtained by Erd\H os and Turan \cite{ErdTur-1941}.
In general $C_2[g]$ sets are those sets $A$ such that each difference $a-a'$ appears at most $g-1$ times and Theorem \ref{thm1} recovers  Corollary 2.1 in \cite{Cill-2010}.

In the opposite direction, Erd\H{o}s and Harzheim  proved that

\lq \lq \emph{if $\alpha$ satisfies $0 < \alpha < \Big (1-\frac 1h\Big )\Big (1-\frac 1g\Big )$,
then for all sufficiently large natural numbers $n$ there exists a subset of $\{1,\cdots, n\}$
 which has at least $n^\alpha$ elements but no $g$ disjoint congruent $h$-element subsets.}\rq \rq \
 
Erd\H os and Harzheim added the restriction that the sets are disjoint because it simplifies the proof. We do the same but obtain a better lower bound. It seems that the same lower bound should hold for $C_h[g]$ sets but we have not found a proof.
\begin{thm} \label{thm2}
There exists a weak-$C_h[g]$ subset $A$ of $[1,n]$ such that
$$|A| \gg n^{\left( 1- \frac{1}{h}\right)\left( 1- \frac{1}{g} \right)\left(1+\frac{1}{hg-1} \right)}
.$$
\end{thm}
It should be noted that for $h$ fixed, Theorem \ref{thm2} gives $|A|\gg n^{1-\frac 1h-\epsilon}$ for $g$ sufficiently large, being the lower bound close to the exponent given in Theorem  \ref{eq:BetterUpperBound-for-F_h(g,n)}. For small values of $g$, Theorem \ref{thm2} is not so strong. Indeed there are algebraic constructions of Sidon sets (which are obviously $C_h[g]$ sets) of size $|A|\sim n^{1/2}$, while the exponent in Theorem \ref{thm2} is not greater than $1/2$ for $h=2$ and for $(h,g)=(3,3)$. Theorem \ref{thm2} gives non trivial lower bounds in any other case.

Theorem \ref{thm1} and  what is known about Sidon sets give the estimates
$$n^{1/2}(1+o(1))< |A|<2^{1/3}n^{2/3}(1+o(1))$$
for a $C_3[3]$ set $A\subset \{1,\dots,n\}$ of maximum cardinality.
Probably there exists an algebraic construction of a $C_3[3]$ set beating the exponent in the lower bound but we have not found it.

We prove Theorem \ref{thm2} using the probabilistic method combined with the deletion technique.
A random finite set, say $S$, may not be weak-$C_h[g]$ but have some \lq \lq blemishes\rq \rq, that is to say,
some arrangements of $g$ congruent disjoint $h$-sets occur inside $S$. Blemishes are more likely the larger is $S$.
We tune the construction so that the size of the random set \emph{moderately} exceeds the threshold for the weak-$C_h[g]$ property.
Then we can prune all the blemishes from $S$ obtaining a true weak-$C_h[g]$ set, while the size is
roughly preserved.
These ideas have appeared before in the literature, see for example \cite[\S 3]{AloSpe-2000}, \cite{SpeTet-1995},  and \cite{Cill-2010a}.

Erd\H{o}s and Harzheim also proved that for any $C_h[g]$ infinite sequence $A$ we have
\begin{equation*} \label{eq:UpperEstimate-C_h[g]-sequence}
\liminf_{x \to \infty} \frac{A(x)}{x^{1-1/h}} = 0.
\end{equation*}
We refine this result as follows.
\begin{thm} \label{thm3}
If $A$ is an infinite $C_h[g]$ sequence then
\begin{equation*}
\liminf_{x \to \infty} A(x)\cdot \frac{(x\log x)^{1/h}}x \ll 1,
\end{equation*}
where the implicit constant depends on $g$ and $h$.
\end{thm}

Theorem \ref{thm3} was proved by Erd\H os \cite{Erd-1957} when $h=g=2$ (infinite Sidon sequences).



\section{Finite $C_h[g]$ sets}
We start recalling a Theorem, which is a consequence of the Jensen's inequality, and that will be used next in the proof of Theorem \ref{thm1}.
\begin{thm}[Overlapping theorem \cite{CillTen-2007}] \label{thm:overlapping-thm}
Let $(\Omega, \mathcal{A}, \mathbb{P})$ be a probability space and let $\{ E_j \}_{j=1}^k$ denote a family of events. Write
\[
\sigma_m := \sum_{1 \le j_1 < \cdots < j_m \le k} \Pr{ E_{j_1} \cap \cdots \cap E_{j_m}}, \quad (m \ge 1).
\]
Then we have
\[
\sigma_m \ge \binom{\sigma_1}{m}=\frac{\sigma_1(\sigma_1-1)\cdots (\sigma_1-(m-1))}{m!}.
\]
\end{thm}

\subsection{Proof of Theorem \ref{thm1}} Let $A$ be a $C_h[g]$ set  and let $B$ be any subset of integers within $[1,n]$ of size at least $h$. 
Let $\Y$ be a random variable with range the positive integers and law
\[
    \Pr{ \Y=m} =
	\begin{cases}
 	\dfrac{1}{|A+B|} & \text{if }  m \in A+B, \\
 	0 & \text{otherwise. }
	\end{cases}
\]
For every $b \in B$ we define the event $E_b = \{ \omega \in \Omega \colon \Y(\omega) \in  A + b \},$ that has probability
$\Pr{E_b}=\sum_{a\in A}\Pr{\Y=a+b}= |A|/|A+B|$. We also write
\[
\sigma_m := \sum_{\{b_{1}, \cdots, b_{m}\} \in \binom{B}{m}} \Pr{ E_{b_1} \cap \cdots \cap E_{b_m} }, \quad (m \ge 1).
\]
In particular
$$
 \sigma_1 = \frac{|A||B|}{|A+B|}.
$$
Let  $b_1 > \cdots > b_h$ be $h$ fixed elements of $B$. We can write
\begin{align*}
 \Pr{E_{b_1} \cap \cdots \cap E_{b_{h}}}  & = \sum_{ \{a_1, \cdots, a_h  \} \in \binom{A}{h} } \Pr{\Y=a_1+b_1=a_2+b_2 \cdots =a_h+b_h} \\
 							& = \sum_{ a_1 + \{0,b_1-b_2, b_1-b_3, \cdots, b_1-b_h \}  \in \binom{A}{h}} \frac{1}{|A+B|},
\end{align*}
the sum extending to all $a_1\in A$ such that
$ a_1 + \{0,b_1-b_2, b_1-b_3, \cdots, b_1-b_h \}\subset A.$ These are congruent $h$-subsets of the  $C_{h}[g]$ set $A$, thus
\begin{align*} 							
							 \Pr{E_{b_1} \cap \cdots \cap E_{b_{h}}}
							& \le \dfrac{g-1}{|A+B|}.
\end{align*}
Now we use Theorem \ref{thm:overlapping-thm}  to obtain
$$
	\binom{|B|}{h} \; \dfrac{g-1}{|A+B|} \ge \sigma_{h} \ge  \frac{\sigma_1 (\sigma_1-1) \cdots  (\sigma_1-h+1)}{h!}\ge \frac{\sigma_1}{h!}(\sigma_1-(h-1))^{h-1},
$$
and so
$$
 \frac{|B|^h}{h!} \; \dfrac{(g-1)}{|A+B|} \ge \frac{|A||B|}{h! \, |A+B|} \left(\frac{|A||B|}{|A+B|}-(h-1) \right)^{h-1},
$$
which implies
\[
    |A|^{h/(h-1)} \le | A + B | \left( (g-1)^{1/(h-1)} + \frac{(h-1)|A|^{1/(h-1)}}{|B|} \right).
\]
If we choose $B = [0,\ell]$, by the last inequality we have
\begin{equation} \label{eq:MejorCota-de-F_h(g,n)_1}
|A|^{h/(h-1)} \le (n+\ell) \left( (g-1)^{1/(h-1)} + \frac{(h-1)|A|^{1/(h-1)}}{\ell+1} \right).
\end{equation}
We first take  $\ell=n$ and use  $|A|\le n$ in the right side, getting $|A|^{h/(h-1)} \ll n \implies |A|^{1/(h-1)}\ll n^{1/h} $.
Inserting this in the second member of \eqref{eq:MejorCota-de-F_h(g,n)_1} we obtain
\[
|A|^{h/(h-1)} \le (g-1)^{1/(h-1)} n + O( \ell) \ + O\left (\frac{ n^{1+1/h}}{\ell+1}\right ) + O( n^{1/h}).
\]
To minimize this last upper bound we choose $\ell \asymp n^{1/2+1/2h}$. Then we can write
\[
|A|^{h/(h-1)} \le (g-1)^{1/(h-1)} n + O \left( n^{1/2+1/2h} \right) = (g-1)^{1/(h-1)} n \left( 1 + O\left( n^{1/2h-1/2} \right)\right),
\]
which yields
\[
|A| \le (g-1)^{1/h}  n^{1-1/h} \left( 1 + O\left( n^{1/2h-1/2} \right)\right)^{1 -1/h} = (g-1)^{1/h} n^{1-1/h} + O\left( n^{1/2 - 1/2h)} \right). \qedhere
\]

\subsection{Proof of Theorem \ref{thm2}} \label{sec:proof-of-thm2}
We say that
$m \in S$ is $(h,g)$-bad (for $S$) if there exist $m_1 < \cdots < m_{g-1}$, with $m_i < m$, and there exist $\ell_1<\ell_2 < \cdots < \ell_{h-1}$ such that the sums $\{m_1,\cdots, m_{g-1},m \} + \{0,\ell_1, \cdots, \ell_{h-1}\} $ are $gh$ distinct elements of $S$.

We define $S_{bad}$ the set of $(h,g)$-bad elements for $S$.
It is clear that for any set $S$, the set $$S_{C_h[g]}=S\setminus{S_{bad}},$$ is weak-$C_h[g]$ with cardinality $|S_{C_h[g]}|=|S|-|S_{bad}|.$

Define $p$ as the number such that $2pn=n^{g+h-1}(2p)^{hg}$. It is straightforward to check that \begin{equation}\label{np}np=\frac 12n^{\left( 1- \frac{1}{h}\right)\left( 1- \frac{1}{g} \right)\left(1+\frac{1}{hg-1} \right)}.\end{equation}
We will prove that except for finitely many $n$ there exist a set $S \subset [1,n]$ such that 
\begin{equation}\label{ine}
|S|\ge\frac{np}2\quad \text{ and }\quad |S_{bad}|\le\frac{np}4.
\end{equation}
Note that for such a set we have $$|S_{C_h[g]}|=|S|-|S_{bad}|>\frac{np}4=\frac{n^{\left( 1- \frac{1}{h}\right)\left( 1- \frac{1}{g} \right)\left(1+\frac{1}{hg-1} \right)} }8,$$ for all sufficiently large $n$ and $A=S_{C_h[g]}$ satisfies the conditions of Theorem \ref{thm2}.

Indeed we will prove that  with probability at least $1/4$, a random set $S$  in $[1,n]$ satisfies \eqref{ine} if each element in $[1,n]$ is independently chosen to be in $S$ with probability $p$.

Next we obtain estimates for the random variables $|S|$ and $|S_{bad}|$.

If $m$ is $(h,g)$-bad then the $gh$ sums $\{m_1,\cdots, m_{g-1}, m \} + \{0,\ell_1, \cdots, \ell_{h-1}\}$ are all distinct elements of $S$ and so
$$\Pr{\{m_1,\dots,m_{g-1},m\}+\{0,\ell_1,\dots,\ell_{h-1}\}\subset S}= p^{gh}.$$ 
Hence
\begin{align*}
&\Pr{m \text{ is $(h,g)$-bad} } \le
\sum_{ \substack{1 \le m_1 < \cdots < m_{g-1} < m \\1 \le \ell_1 < \cdots < \ell_{h-1} \le n } }  p^{gh}\le \binom m{g-1}\binom n{h-1}p^{gh}<n^{g+h-2}p^{gh}
\end{align*}
\[
\implies \Esp{|S_{bad}|} \le \sum_{1 \le m \le n} \Pr{m \text{ is $g$-bad} } \le n^{g+h-1}p^{gh} .
\] 
On the one hand by Markov's inequality we have
\begin{align}\label{Markov}
\Pr{|S_{bad}|>\frac{np}4}	&=   \Pr{|S_{bad}|>\frac{n^{g+h-1}(2p)^{gh}}8 }\\
										&=   \Pr{|S_{bad}|>2^{gh-3}n^{g+h-1}p^{gh}   }\nonumber\\
										&\le \Pr{|S_{bad}|>2\Esp{|S_{bad}|}    }<1/2.\nonumber
\end{align}
On the other hand, using that $\Esp{|S|}=np$ and $\text{Var}(|S|)=np(1-p)$ and applying Chebychev's inequality we have
\begin{eqnarray}\label{Che}
\Pr{|S|<\frac{np}2}&=&\Pr{|S|<\frac{\Esp{|S|}}2}<\Pr{|S-\Esp{|S|}|>\frac{\Esp{|S|}}2}\\ &<&\frac{4\text{Var}(|S|)}{(\Esp{|S|})^2}=\frac{4np(1-p)}{(pn)^2}<\frac{4}{pn}<\frac 14,\nonumber
\end{eqnarray} except for finitely many $n$.
By \eqref{Markov} and \eqref{Che} we have
$$ \Pr{|S|\ge np/2 \text{ and } |S_{bad}|\le np/4}\ge 1-\left (1/2+1/4\right )\ge 1/4, $$
as we wanted.
\section{Infinite $C_h[g]$ sequences}

\subsection{Proof of Theorem \ref{thm3}}

For a positive integer $N$, let $[0,N^2]$ denote all the positive integers less or equal to $N^2.$
We divide $[0,N^2]$ into equally sized intervals $$I_\nu := [(\nu-1)N, \nu N], \; \nu=1, \cdots, N.$$
Let $\mathcal{C}$ denote the collection of all $h$-subsets of $[0,N^2]$ that are included in one of the intervals $I_\nu$:
\[
 \mathcal{C} := \left \lbrace C \in \binom{[0,N^2]}{h} \colon C \subset I_\nu \text{ for some } \nu \right \rbrace.
\]
We say that the sets in the collection $\mathcal{C}$ are \lq\lq small\rq\rq \ as their diameter is at most $N$.
We classify the elements of $\mathcal{C}$ so that each class groups all the sets that are pairwise congruent.
Each class $\alpha $ contains a set $C_\alpha$ that contains $0$,
and the remaing $h-1$ elements of $C_\alpha$ can be chosen in $\binom{N-1}{h-1}$ different ways;
each of the choices determines a class different from the others. Then the number of classes is
\[
  \binom{N-1}{h-1}.
\]
Let $A_\nu$ denote the size of $A \, \cap \, I_\nu$, we have $A_\nu = A(\nu N) - A((\nu-1)N)$, where $A(x) := |\{ a \in A \colon a \le x \}|$ is the counting function of the sequence.

One the one hand as $A$ is a $C_h[g]$ sequence then in every class of $\mathcal{C}$ there are at most $g-1$ subsets of $A$. 
Hence we have the following upper bound for the total number of \lq\lq small\rq\rq \ subsets of $A$ that belong to $\mathcal{C}$
\begin{equation*} \label{eq:liminf-ll-1-1}
\sum_{\nu=1}^{N} \binom{A_\nu}{h} \le \binom{N-1}{h-1} (g-1) \ll N^{h-1} \quad (N \to \infty),
\end{equation*}
Now we prove by induction in $h$ that
\begin{equation} \label{eq:liminf-ll-1-3}
\sum_{\nu=1}^{N} A_\nu^h \ll N^{h-1} \quad (N \to \infty).
\end{equation}
For $h=2$ we know by Theorem \ref{eq:BetterUpperBound-for-F_h(g,n)} that $A(N^2)\ll N$, so
$$\sum_{\nu=1}^{N} A_\nu^2 = 2 \sum_{\nu=1}^{N} \binom{A_\nu}{2} + \sum_{\nu=1}^{N} A_\nu \ll N + A(N^2) \ll N.$$
If \eqref{eq:liminf-ll-1-3} holds for all exponents up to $h-1$, then
\begin{equation*}
\sum_{\nu=1}^{N} A_\nu^h = h! \sum_{\nu=1}^{N} \binom{A_\nu}{h}+ O\left( \sum_{\nu=1}^{N} A_\nu^{h-1} \right) \ll N^{h-1} +N^{h-2}, \qquad (N \to \infty),
\end{equation*}
thus it also holds for $h$. Using \eqref{eq:liminf-ll-1-3} and  H\"{o}lder inequality we can write
\begin{align}
\sum_{\nu=1}^N A_\nu \left(\frac{1}{\nu}\right)^{1-1/h} & \le \left( \sum_{\nu=1}^N A_\nu^h \right)^{1/h} \, \left( \sum_{\nu=1}^N \frac{1}{\nu} \right)^{1-1/h} \notag \\
		&\ll \left( N \log N \right)^{1-1/h}, \qquad (N \to \infty) \label{eq:liminf-ll-1-4}.
\end{align}

\medskip
On the other hand as $\sum_{\nu \le t} A_\nu=A(tN)$ and summing by parts
\[
\sum_{\nu=1}^N A_\nu \left(\frac{1}{\nu}\right)^{1-1/h} \gg \int_1^N \frac{A(tN)}{t^{2-1/h}} \ud t.
\]
Let us write $$\tau(m):=\inf_{n \ge m} \frac{A(n)(\log n)^{1/h}}{n^{1-1/h}}.$$
For $N \ge m$ we have
\[
\int_1^N \frac{A(tN)}{t^{2-1/h}} \ud t \gg \frac{\tau(m)N^{1-1/h}}{(\log N)^{1/h}} \int_1^N \frac{1}{t} \ud t \gg \tau(m) (N \log N)^{1-1/h},
\]
and so
\[
\sum_{\nu=1}^N A_\nu \left(\frac{1}{\nu}\right)^{1-1/h} \gg \tau(m) (N \log N)^{1-1/h}.
\]
Inserting $\eqref{eq:liminf-ll-1-4}$ we have $\lim_{m \to \infty} \tau(m) \ll 1,$
that is what we wanted to prove.

\section{Open problems}

Theorem \ref{thm1} and  what is known about Sidon sets give the estimates
$$n^{1/2}(1+o(1))< |A|<2^{1/3}n^{2/3}(1+o(1))$$
for a $C_3[3]$ set $A\subset \{1,\dots,n\}$ of maximum cardinality.

\textbf{Problem 1:} Improve the bounds above.

Probably there exists an algebraic construction of a $C_3[3]$ set beating the exponent in the lower bound but we have not found it.

The analogous problem  in two dimensions can be nicely illustrated.

\textbf{Problem 2:} What is the largest size of a set $A\in [1,n]\times [1,n]$ avoiding three translated triangles?

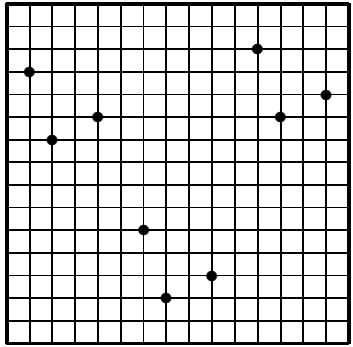
\setlength{\unitlength}{3mm}
\begin{picture}(0,0.3)
   { \linethickness{.3mm}
  \put(10, 0){\line(1, 0){15}}
  \put(10, -15){\line(1, 0){15}}
  \put(10, 0){\line(0, -1){15}}
  \put(25, 0){\line(0, -1){15}}}
  \linethickness{.1mm}
   \put(10, -1){\line(1, 0){15}}
  \put(10, -2){\line(1, 0){15}}
  \put(10, -3){\line(1, 0){15}}
  \put(10, -4){\line(1, 0){15}}
  \put(10, -5){\line(1, 0){15}}
  \put(10, -6){\line(1, 0){15}}
  \put(10, -7){\line(1, 0){15}}
  \put(10, -8){\line(1, 0){15}}
   \put(10, -9){\line(1, 0){15}}
  \put(10, -10){\line(1, 0){15}}
  \put(10, -11){\line(1, 0){15}}
  \put(10, -12){\line(1, 0){15}}
  \put(10, -13){\line(1, 0){15}}  
   \put(10, -14){\line(1, 0){15}} 
    \put(11, 0){\line(0, -1){15}}
     \put(12, 0){\line(0, -1){15}}
      \put(13, 0){\line(0, -1){15}}
     \put(14, 0){\line(0, -1){15}}
     \put(15, 0){\line(0, -1){15}}
     \put(16, 0){\line(0, -1){15}}
      \put(17, 0){\line(0, -1){15}}
     \put(18, 0){\line(0, -1){15}}
     \put(19, 0){\line(0, -1){15}}
     \put(20, 0){\line(0, -1){15}}
      \put(21, 0){\line(0, -1){15}}
     \put(22, 0){\line(0, -1){15}}
     \put(23, 0){\line(0, -1){15}}
     \put(24, 0){\line(0, -1){15}}
      \put(25, 0){\line(0, -1){15}}        
 {\put(11,-3){\circle*{0.5}}
   \put(12,-6){\circle*{0.5}}
    \put(14,-5){\circle*{0.5}}
   \put(16,-10){\circle*{0.5}}
   \put(17,-13){\circle*{0.5}}
    \put(19,-12){\circle*{0.5}}
      \put(21,-2){\circle*{0.5}}
   \put(22,-5){\circle*{0.5}}
    \put(24,-4){\circle*{0.5}}}
 \end{picture}
 
 \vspace{5cm}

Our method would only give $n(1+o(1))\le|A|\ll n^{4/3}$.
 
\textbf{Problem 3:} Remove the condition \emph{weak} in Theorem \ref{thm2}.
 
\textbf{Problem 4:} Construct infinite dense $C_h[g]$ sequences. At least with counting function $A(x)\gg x^{\left( 1- \frac{1}{h}\right)\left( 1- \frac{1}{g} \right)\left(1+\frac{1}{hg-1} \right)+o(1)}$.

We have found technical difficulties to deal with problems 3 and 4. Erd\H os and Harzheim probably also found these difficulties, which would explain why they added the weak condition in his lower bound for the finite case and did not include the infinite case in their study.

\end{document}